\documentclass[12pt]{article}
\usepackage{amsmath}
\usepackage{amssymb} 
\usepackage{latexsym}
\usepackage{theorem}
\usepackage{euscript}
\newcommand{\msf}[1]{{\mathsf {#1}}}

\newcommand{\mcal}[1]{{\mathcal {#1}}} 

\newcommand{\hh}{\ensuremath{\mathfrak h}}
\renewcommand{\gg}{\ensuremath{\mathfrak g}}
\renewcommand{\ss}{\ensuremath{\mathfrak s}}
\newcommand{\co}{\colon\thinspace} 

\newtheorem{theorem}{Theorem}  

\newtheorem{lemma}[theorem]{Lemma}
\newtheorem{corollary}[theorem]{Corollary}

{\theoremstyle{break}	
{\theorembodyfont{\rmfamily} 
				\newtheorem{example}[theorem]{Example}
				}}
\def\qed{\ifhmode\unskip\nobreak\fi\ifmmode\ifinner\else\hskip5pt\fi\fi
 \hfill\hbox{\hskip5pt\vrule width4pt height6pt depth1.5pt\hskip1pt}}
\newenvironment{proof}[1]{\smallskip \noindent {\bf #1}}{\qed\smallskip}

\newcommand{\RR}{\ensuremath{{\mathbb R}}}     
\newcommand{\R}[1]{\ensuremath{{\mathbb R}^{#1}}} 
\newcommand{\NN}{\ensuremath{{\mathbb N}}}     

\renewcommand{\P}[1]{\ensuremath{{\bf P}^{#1}}} 

\newcommand{\CC}{\ensuremath{{\mathbb C}}}     

      \def\mylabel#1{\label{#1}}

\begin{document}
\author{\sc Morris W. Hirsch\\
 University of California at Berkeley\\
 University of Wisconsin  at Madison \and
\sc Joel W. Robbin\\ University of Wisconsin  at Madison}
\title{ON EIGENVECTORS OF NILPOTENT LIE ALGEBRAS OF LINEAR OPERATORS}
\maketitle

\begin{abstract}
   We give a condition ensuring that the operators in a nilpotent Lie
algebra of linear operators on a finite dimensional vector space have
a common eigenvector.
\end{abstract}
\hspace*{2.3in}{\sc Introduction}

\smallskip
Throughout this paper $V$ is a vector space of positive dimension over
a field $\msf f$ and $\gg$ is a nilpotent Lie algebra over $\msf f$ of
linear operators on  $V$.  An element $u\in V$  is an {\em eigenvector
for $S\subset  \gg$} if  $u$ is an  eigenvector for every  operator in
$S$.  If $V$ has a  basis $(e_1,\dots, e_n)$ representing each element
of $\gg$ by  an upper triangular matrix, then  $e_1$ is an eigenvector
for $\gg$.  Such a basis  exists when $\msf f$ is algebraically closed
and $\gg$ is solvable (Lie's  Theorem), and also when every element of
$\gg$  is a  nilpotent operator  (Engel's Theorem).   Our  results are
further conditions guaranteeing existence eigenvectors.

 The minimal and characteristic polynomials of a linear operator $A$
on $V$ are denoted respectively by $\pi_A, \mu_A \in\msf f[t]$ = the
ring of polynomials over $\msf f$.  
written $\#S$.

Let $\msf k$ be a Galois extension field of $\msf f$ of degree
$d:=[\msf k:\msf f]$,  and  define ${\bf 
M}\subset\NN$ to be the additive monoid generated by zero and 
the prime divisors $d$.

Consider the conditions:
\begin{description}
\item[(C1)] $\mu_A$ splits in $\msf k$ for every $A\in\gg$

\item[(C2)] $\dim V\notin {\bf M}$
 
\end{description}

Our main result is:
\begin{theorem}		\mylabel{th:A}
If {\bf (C1)} and {\bf (C2)} hold then $\gg$ has an eigenvector. 
\end{theorem}
The proof is preceded by some applications.

\medskip
When (C1) holds, Theorem  \ref{th:A} shows that  there is an
eigenvector in every  invariant subspace whose dimension is not in
$\bf M$.  This is exploited to yield the following two results:
\begin{corollary}		\mylabel{th:corB}
If a nilpotent Lie algebra  of linear operators on $\R n$ does
not have an eigenvector,  every nontrivial invariant subspace has odd
dimension.
\end{corollary}
\begin{proof} When $\msf f$ is the real field $\RR$  and $\msf k$ is
the complex field $\CC$,  $\ \bf M$ consists of the positive even
integers.
\end{proof}

\begin{corollary}		\mylabel{th:D}
Let {\bf (C1)} hold.  Assume $\gg$ preserves a direct sum
decomposition $V=\oplus_i W_i$, and let $D\subset \NN$ denote the set
of dimensions of the subspaces $W_i$.
\begin{description}

\item[(i)] If $\gg$ does not have an eigenvector then  $D\subset\bf M$.

\item[(ii)] If $V'\subset V$ is a maximal subspace spanned by
eigenvectors of $\gg$ then $\dim (V')\ge \#\{D\setminus \bf M\}$.

\end{description}
\end{corollary}
\begin{proof}
Assertion (i) follows from Theorem \ref{th:A}.  To prove (ii) order
the $W_i$ so that $W_1,\dots,W_m$ are the only summands whose
dimensions are not in $\bf M$.  For each $j\in\{1,\dots, m\}$ we
choose an eigenvector $e_j\in W_j$ by Theorem \ref{th:A}.  The $e_j$
are linearly independent and belong to $V'$ by maximality of $V'$,
whence (ii).
\end{proof}

\begin{example}		\mylabel{th:commuting}
Assume $n\notin \bf M$ and let $\alpha\in \msf f[t]$ be a monic
polynomial that splits in $\msf k[t]$.  Denote by $\mcal A (\alpha)$
the set of $n\times n$ matrices $T$ over $\msf f$ such that $\alpha
(T)=0$.  Then every pairwise commuting family $\mcal T\subset \mcal A
(\alpha)$ has an eigenvector in $\msf f^{n}$.  This follows from
Theorem \ref{th:A} applied to the Lie algebra $\gg$ of linear
operators on $\msf f^n$ generated by $\mcal T$.  Being abelian, $\gg$
can be triangularized over $\msf k$, hence (C1) holds.
\end{example}

\begin{example}		\mylabel{th:ex2}
The assumption that $n \in {\bf M}$ is essential to Theorem
\ref{th:A}.  For instance, take $\msf f=\RR$, $\msf k=\CC$, $V=\R 2$. 
 The abelian  Lie algebra of $2\times 2$ of real skew symmetric matrices.
does not have an eigenvector in $\R 2$.
\end{example}

\begin{example}		\mylabel{th:ex4}
The hypothesis of Theorem \ref{th:A} cannot be weakened to $\gg$ being
merely solvable.  For a
 counter\-example with $\msf f=\RR, \msf k=\CC$, take $\gg$ to be the
 solvable 3-dimensional real Lie algebra with basis $(X, U, V)$ such
 that $[X,U]=-V, \; [X,V]=U,\; [U, V]=0$.
\end{example}
 
A Lie algebra $\ss$ over $\msf f$ is {\em supersolvable} if the
spectrum of the linear map ${\rm ad}\ A\co \ss\to\ss$ lies in $\msf f$
for all $A\in\ss$.  If $\ss$ is not supersolvable it need not have an
eigenvector, as is shown by Example \ref{th:ex4}.  We don't know if
Theorem \ref{th:A} extends to supersolvable Lie algebras, except for
the following special case:

\begin{theorem}		\mylabel{th:super}
A supersolvable Lie algebra $\ss$ of linear transformations of
$\R 3$ has an eigenvector.
\end{theorem}
\begin{proof}
 Lacking an algebraic proof, we use a dynamical argument.  Let
$G\subset GL(3,\RR)$ be the connected Lie subgroup having Lie algebra
$\ss$.  The natural action of $G$ on the projective plane $\P 2$ of
lines in $\R 3$ through the origin fixes some $L\in \P 2$.  This
follows from supersolvability because $\dim (\P 2) =2$, the action on
$\P 2$ is effective and analytic, and the Euler characteristic of $\P
2$ is nonzero (Hirsch \& Weinstein \cite{HW00}).  The nonzero points
of $L$ are eigenvectors for $\ss$.
\end{proof}

\medskip

\hspace*{1.7in}{\sc Proof of Theorem \ref{th:A}}

\smallskip
We rely on Jacobson's {\em Primary Decomposition Theorem} \cite[II.4,
Theorem 5] {Jacobson62}.  This states that $V$ has a $\gg$-invariant
direct sum decompsition $\oplus V_i$ where each {\em primary
component} $V_i$ has the following property: For each $A\in\gg$ the
minimal polynomial of $A|V_i$ is a prime power in $\msf f [t]$.

Condition (C2) implies the dimension of some primary component is
$\notin \bf M$.  To prove Theorem \ref{th:A} it therefore suffices to
apply the following result to such a primary component:

\begin{theorem}		\mylabel{th:A'}
Assume {\bf (C1)} and {\bf (C2)}.  If $\pi_A$ is a prime power in
$\msf f[t]$ for each $A\in \gg$ then the following hold:
\begin{description}

\item[(a)] $\pi_A(t)=(t-r_A)^n,\: r_A\in\msf f$

\item [(b)] there is a basis putting  $\gg$ in triangular form
\end{description}
 \end{theorem}

Assertion (a) is equivalent to $\pi_A$ having a root $r_A\in\msf f$.
Therefore (a) follows from:

\begin{lemma}		\mylabel{th:C}
Let $\alpha\in \msf f[t]$ be a polynomial of degree $n$ that splits in
$\msf k[t]$.  If $n\notin \bf M$ then $\alpha$ has a root in $\msf f$,
and the sum of the multiplicities of such roots is $\notin \bf M$.
\end{lemma}
\begin{proof} 
Let $R\subset\msf k$ denote the set of roots of $\pi$, and
$R_j\subset R$  the set of roots of multiplicity $j$. 
The Galois group $\Gamma$ has order $[\msf k:\msf f]$ and acts on $R$ by
permutations.  The cardinality of each orbit divides $[\msf k:\msf
f]$, and $R\cap \msf f$ is the set of fixed points of this action.
Each $R_j$ is a union of orbits, as is $R_j  \verb=\=  \msf f$.
It follows that $\#(R_j \ \verb=\= \ \msf f)\in {\bf M}$.

Let $k\le n$ denote the sum of the multiplicities of the roots that
are not in $\msf f$.  Then 
\[
  k=\sum_{j=2}^n  j\cdot\#(R_j \ \verb=\= \ \msf f)
\]
Therefore $k\in \bf M$ because ${\bf M}$ is closed under addition. 
By hypothesis $n\notin {\bf M}$, hence $n-k\notin {\bf M}$ and
$n-k >0$.  As $n-k$
is the sum of the multiplicities of the roots in $\msf f$, the
conclusion follows.
\end{proof}

Now that  (a) of Theorem \ref{th:A'} is proved, assertion (b) is a
consequence of the following result: 

\begin{lemma}		\mylabel{th:lem1}
Let $\hh$ be a nilpotent Lie algebra of linear operators on $V$. 
Assume that for all $A\in\hh$ there exists
$r_A\in\msf f$ such that 
$  \pi_A (t)= (t-r_A)^n$. 
Then $V$ has a basis putting $\hh$ in triangular form. 
\end{lemma}
\begin{proof}
Every $A\in\hh$ can be written uniquely as $r_AI + N_A$ with $N_A$
nilpotent and $I$ the identity map of $V$.  It is easy to see that the
set comprising the $N_A$ is closed under commutator brackets.
Therefore $V$ has a basis triangularizing all the $N_A$ (Jacobson
\cite[II.2, Theorem 1$'$] {Jacobson62}), and such a basis
triangularizes $\hh$.
\end{proof}

This completes the proof of Theorem \ref{th:A}.


\begin{thebibliography}{9}

\bibitem{HW00} M. Hirsch \& A. Weinstein, {\em Fixed points of
analytic actions of supersoluble Lie groups on compact surfaces},
Ergod. Th. Dyn. Sys.  {\bf 21} (2001), no. 6, 1783--1787.  See also 
{\tt http://front.math.ucdavis.edu/math.DS/0002013}


\bibitem{Jacobson62} N. Jacobson, ``Lie Algebras'', Interscience Tracts
  in Pure Mathematics No. 10.  John Wiley, New York (1962)

\end{thebibliography}
\end{document}